\documentclass[12pt,reqno]{amsart}
\usepackage{amsxtra}
\usepackage{amssymb}
\usepackage{amscd}
\usepackage{mathrsfs}
\newtheorem{theorem}{Theorem}[section]
\newtheorem{lemma}[theorem]{Lemma}

\theoremstyle{definition}

\theoremstyle{remark}

\newtheorem{remark}{Remark}[section]
\newcommand*{\ptens}[1]{\mathop{\widehat\otimes}_{#1}}
\newcommand*{\Ptens}{\mathop{\widehat\otimes}}

\newcommand*{\lmod}{\mbox{-}\!\mathop{\mathbf{mod}}}
\newcommand*{\rmod}{\mathop{\mathbf{mod}}\!\mbox{-}}
\newcommand*{\bimod}{\mbox{-}\!\mathop{\mathbf{mod}}\!\mbox{-}}
\newcommand*{\id}{\mathbf 1}
\newcommand*{\op}{\mathrm{op}}

\newcommand*{\h}{\mathbf h}
\DeclareMathOperator{\Ker}{Ker}
\DeclareMathOperator{\Tr}{Tr}

\newcommand{\lriso}{\stackrel{\textstyle\sim}{\smash\longrightarrow%
\vphantom{\scriptscriptstyle{_1}}}}

\newcommand*{\CC}{\mathbb C}

\newcommand*{\A}{\mathfrak A}
\newcommand*{\nN}{\mathscr N}
\newcommand*{\sS}{\mathcal S}
\newcommand*{\B}{\mathscr B}
\renewcommand*{\H}{\mathscr H}
\newcommand*{\rar}{\rightarrow}

\newcommand*{\xra}{\xrightarrow}
\begin{document}
\title[Convolution algebras
of nuclear operators]{Biprojectivity and biflatness for convolution algebras
of nuclear operators}
\author{A. Yu. Pirkovskii}
\thanks{Partially supported by the RFBR grants 02-01-00928
and 01-01-00490}
\date{}
\begin{abstract}
For a locally compact group $G$, the convolution product on
the space $\nN(L^p(G))$ of nuclear operators was defined by Neufang
\cite{Neuf_PhD}.
We study homological properties of the convolution algebra
$\nN(L^p(G))$ and relate them with some properties of the group $G$,
such as compactness, finiteness, discreteness, and amenability.
\end{abstract} 
\maketitle
\section{Introduction}
Let $G$ be a locally compact group and let $H=L^2(G)$.
In \cite{Neuf_PhD}, M.~Neufang defined a new product on the
space $\nN(H)$ of nuclear operators on $H$ making it
into a Banach algebra. While the usual product on $\nN(H)$ can be
viewed as a noncommutative version of the pointwise product on $\ell^1$, this
new product is an analogue of the convolution product on $L^1(G)$.
The resulting Banach algebra $\sS_1(G)=(\nN(H),*)$
shares many properties with $L^1(G)$ \cite{Neuf_PhD}. In particular,
the group $G$ can be completely reconstructed from $\sS_1(G)$
(compare with the classical result of Wendel \cite{Wendel} about $L^1(G)$).
Some important theorems of harmonic analysis (e.g., a theorem of Hewitt
and Ross \cite[35.13]{HRII} characterizing multipliers on $L^\infty(G)$ for a compact
group $G$) also have their counterparts for $\sS_1(G)$.
Further, $\sS_1(G)$ (like $L^1(G)$) has a right identity if and only if $G$ is discrete,
always has a right b.a.i., is a right ideal in its bidual if and only if $G$ is compact.
A closed subspace of $\sS_1(G)$ is a right ideal if and only if it is invariant
w.r.t. a certain natural action of $G$. Other examples showing that
$\sS_1(G)$ behaves in much the same way as $L^1(G)$ can be found in
\cite{Neuf_PhD}.   

The aim of this paper is to study some homological properties of $\sS_1(G)$,
specifically, {\em biprojectivity} and {\em biflatness}. 
(For a detailed exposition of the homology theory for Banach algebras
we refer to \cite{X1};
some facts can also be found in \cite{X2}, \cite{Dales}, and \cite{Runde}).
Recall that $L^1(G)$ is
biprojective if and only if $G$ is compact \cite{X_calc_est,X_flat},
and is biflat if and only if $G$ is amenable
\cite{Jhnsn_CBA,X_flat}. Thus it is natural
to ask whether similar results hold for $\sS_1(G)$. We show
(Theorem \ref{thm:biflat}) that the latter result concerning the biflatness
of $L^1(G)$ is also true for $\sS_1(G)$. On the other hand, it turns out
(Theorem \ref{thm:bipr}) that $\sS_1(G)$ is biprojective if and only if $G$
is finite. We also show that properties of $G$ such as discreteness and
compactness are equivalent to projectivity of certain $\sS_1(G)$-modules.

\begin{remark}
Perhaps it is appropriate to note that $\sS_1(G)$ is never amenable
(except for the trivial case $G=\{ e\}$) because it has a nontrivial right annihilator
(see \cite{Neuf_PhD} or Section~2 below).
\end{remark}

\section{Preliminaries}

Let $G$ be a locally compact group equipped with a left Haar measure,
and let $1 < p <\infty$.
Given $t\in G$, denote by $L_t\colon L^p(G)\to L^p(G)$ the left
translation operator defined by $(L_t f)(s)=f(ts)$.
For a function $f$ on $G$ we set, as usual, $\check f(t)=f(t^{-1})$ and
$\tilde f(t)=\Delta(t^{-1}) f(t^{-1})$, where $\Delta$ is the modular function on $G$.
For each $h\in L^\infty(G)$ we denote by $M_h\colon L^p(G)\to L^p(G)$
the multiplication operator $f\mapsto hf$. The trace duality
between the space of bounded operators,
$\B(L^p(G))$, and the space of nuclear operators, $\nN(L^p(G))$, will be denoted by
the brackets $\langle \:\cdot\: , \:\cdot\: \rangle$.

The {\em convolution product} $*$ on $\nN(L^p(G))$ introduced by Neufang
\cite{Neuf_PhD} is defined as follows. First consider the bilinear map
\begin{gather}
\B(L^p(G))\times\nN(L^p(G))\to L^\infty(G),\notag\\
\label{B_N_L}
(T,\rho)\mapsto \Bigl( t\mapsto \langle\rho, L_t T L_{t^{-1}}\rangle\Bigr)\, .
\end{gather}
Next consider the representation
\begin{equation}
\label{L_B}
L^\infty(G)\to\B(L^p(G)),\qquad h\mapsto M_h\, .
\end{equation}
Composing \eqref{B_N_L} and \eqref{L_B}, we obtain a bilinear map
\begin{equation*}
\B(L^p(G))\times\nN(L^p(G))\xra{\odot} \B(L^p(G)),\qquad 
(T,\rho)\mapsto T\odot\rho\, .
\end{equation*}
For every $\rho\in\nN(L^p(G))$ the map $T\mapsto T\odot\rho$ is
weak$^*$ continuous \cite{Neuf_PhD}.
Therefore we have a well-defined bilinear map
\begin{gather*}
\nN(L^p(G))\times \nN(L^p(G))\xra{*} \nN(L^p(G)),\\
\langle T,\rho * \tau\rangle = \langle T\odot\rho,\tau\rangle
\qquad\text{for each}\quad T\in\B(L^p(G))\, .
\end{gather*}
Neufang (\cite{Neuf_PhD}, Satz~5.2.1 and Prop.~5.4.1)
proved that $(\nN(L^p(G)),*)$ is an associative
Banach algebra with a right b.a.i.
We shall denote this algebra by $\nN^p(G)$. 
(Note that the algebra $\nN^2(G)$ is denoted by $\sS_1(G)$ in \cite{Neuf_PhD}, 
and the dual module $\sS_1^*(G)=\B(L^2(G))$ is denoted by $\sS_\infty(G)$
in order to emphasize a relation with the Schatten classes.)

The algebra $\nN^p(G)$ can be considered as an extension of $L^1(G)$
in the following way. Consider the map $\iota^1\colon L^\infty(G)\to\B(L^p(G))$
defined by the rule $\iota^1(h)=M_{\check h}$. This map is continuous w.r.t. the
weak$^*$ topologies determined by the dualities
$\bigl\langle L^\infty(G),L^1(G)\bigr\rangle$ and
$\bigl\langle \B(L^p(G)),\nN(L^p(G))\bigr\rangle$.
Hence there exists the predual map $\sigma\colon\nN(L^p(G))\to L^1(G)$.
Neufang (\cite{Neuf_PhD}, Satz~5.3.1) proved that $\sigma$ is a Banach algebra
homomorphism from $\nN^p(G)$ onto $L^1(G)$. Furthermore, 
\begin{equation}
\label{Ker_sigma}
\Ker\sigma=\{\rho\in\nN^p(G) \, ; \; \langle \rho,M_h\rangle = 0
\quad\forall\; h\in L^\infty(G)\}\, .
\end{equation}
Therefore (see \cite{Neuf_PhD}, Satz~5.3.4) we have an extension
\begin{equation}
\label{neuf_ext}
0 \rar I \rar \nN^p(G) \xra{\sigma} L^1(G) \rar 0
\end{equation}
of Banach algebras. Note also that the definition of the product in $\nN^p(G)$
together with \eqref{Ker_sigma} implies that $\nN^p(G) I = 0$.

We shall need a more explicit description of $\sigma$.
Take $q\in (1,+\infty)$ such that $\frac{1}{p}+\frac{1}{q}=1$, and
recall that there exists an isometric isomorphism
\begin{equation}
\label{tens_nucl}
L^p(G)\Ptens L^q(G) \lriso \nN(L^p(G)),\qquad
f\otimes g\mapsto \bigl(h\mapsto \langle g,h\rangle f\bigr)\, . 
\end{equation}
Here the brackets $\langle \:\cdot\: , \:\cdot\: \rangle$ denote the
usual $L^p$-$L^q$ duality. Identifying an elementary tensor
$f\otimes g\in L^p(G)\Ptens L^q(G)$ with the corresponding rank-one operator,
we see that
\begin{multline*}
\langle\sigma(f\otimes g),h\rangle 
= \langle f\otimes g,M_{\check h}\rangle
= \langle M_{\check h}(f),g\rangle
= \langle \check h f, g\rangle\\
= \int_G h(t^{-1}) f(t) g(t) \, dt
= \int_G \Delta(t^{-1}) h(t) f(t^{-1}) g(t^{-1}) \, dt
= \langle (fg)^\sim, h\rangle
\end{multline*}
for every $h\in L^\infty(G)$. Therefore,
\begin{equation}
\label{sigma}
\sigma(f\otimes g)=(fg)^\sim \qquad (f\in L^p(G),\; g\in L^q(G))\, .
\end{equation}

Consider the algebra homomorphism $\varepsilon\colon L^1(G)\to\CC$
given by $\varepsilon(f)=\int_G f\, d\mu$. It is clear from \eqref{sigma} that
$\varepsilon\sigma=\Tr$. Thus $\CC$ can be viewed as a $L^1(G)$-module
via $\varepsilon$ and as a $\nN^p(G)$-module via $\Tr$. We shall denote
these modules by $\CC_\varepsilon$ and $\CC_{\Tr}$, respectively.

Recall some notation and some definitions from the homology theory
of Banach algebras (for details, see \cite{X1,X2}).
Let $A$ be a Banach algebra. The category of left (resp. right) Banach
$A$-modules is denoted by $A\lmod$ (resp. $\rmod A$).
If $B$ is another Banach algebra, then $A\bimod B$ stands for the category
of Banach $A$-$B$-bimodules. 
Spaces of morphisms in the above categories
are denoted by ${_A}\h(X,Y)$, $\h_A(X,Y)$, and ${_A}\h_B(X,Y)$, respectively.
The space of continuous linear operators between Banach spaces $X$ and $Y$
is denoted by $\B(X,Y)$.
For each left Banach $A$-module $X$ denote by $A\cdot X\subset X$ the closed
linear span of $\{ a\cdot x\, ;\, a\in A,\; x\in X\}$.
If $A\cdot X=X$, then $X$ is said to be {\em essential}.

A sequence $X_\bullet=(0\to X_1 \to X_2 \to X_3 \to 0)$
of left Banach $A$-modules is called
{\em admissible} (resp. {\em weakly admissible}) if it splits as a
sequence of Banach spaces (resp. if the dual sequence
$X_\bullet^*=(0\to X_3^*\to X_2^* \to X_1^*\to 0)$
splits as a sequence of Banach spaces).
A left Banach $A$-module $Y$ is said to be {\em projective} (resp. {\em injective})
if for each admissible sequence $X_\bullet$ in $A\lmod$
the induced sequence ${_A}\h(Y,X_\bullet)$ (resp. ${_A}\h(X_\bullet,Y)$)
is exact. A left Banach $A$-module $Y$ is called
{\em flat} if for each admissible sequence
$X_\bullet$ in $\rmod A$ the induced sequence
$X_\bullet\ptens{A}Y$ is exact. 
Recall that each projective module is flat, and that
$Y\in A\lmod$ is flat iff the dual module, $Y^*$, is injective in $\rmod A$
\cite[7.1]{X2}, \cite[VII.1]{X1}.
If the {\em canonical morphism} $\pi\colon A\Ptens Y\to Y,\;
a\otimes y\mapsto a\cdot y$ is a retraction in $A\lmod$ (i.e., if there
exists an $A$-module morphism $\rho\colon Y\to A\Ptens Y$
such that $\pi\rho=\id_Y$), then $Y$ is projective.
The converse is true provided $Y$ is essential \cite[IV.I]{X1}.

\begin{remark}
\label{rem:pullback}
Let $A\to B$ be a Banach algebra homomorphism with dense range.
Assume $Y\in B\lmod$ is projective (resp. injective, resp. flat) in $A\lmod$.
Then $Y$ is projective (resp. injective, resp. flat) in $B\lmod$.
To see this, it suffices to observe that
${_A}\h(X,Y)={_B}\h(X,Y)$ for each $X\in B\lmod$, and that
$X\ptens{A}Y\cong X\ptens{B}Y$ for each $X\in\rmod B$.
See also \cite[IV.I]{X1}.
\end{remark} 
 
Projective and flat right $A$-modules and $A$-bimodules
are defined similarly. 

A Banach algebra $A$ is called {\em biprojective} (resp. {\em biflat})
if $A$ is a projective (resp. flat) Banach $A$-bimodule.
Recall that $A$ is biprojective (resp. biflat) iff the product map
$\pi_A\colon A\Ptens A\to A,\; a\otimes b\mapsto ab$ is a retraction in $A\bimod A$
(resp. iff the dual map $\pi_A^*\colon A^*\to (A\Ptens A)^*$ is a coretraction
in $A\bimod A$); see \cite[IV.5 and VII.2]{X1}.

A Banach algebra $A$ is said to be {\em contractible} (resp. {\em amenable})
if the first Hochschild cohomology group, $\H^1(A,X)$, is trivial for each
(resp for each dual) Banach $A$-bimodule $X$.
Recall that $A$ is contractible iff it is biprojective and unital, and
is amenable iff it is biflat and has a b.a.i. \cite[7.1]{X2}.

Following Selivanov \cite{Sel_super}, we say that a Banach algebra $A$
is {\em superbiprojective} (resp. {\em superbiflat}) if it is biprojective
(resp. biflat) and $\H^2(A,X)=0$ for each (resp. for each dual)
Banach $A$-bimodule $X$. Selivanov proved that $A$ is superbiflat
iff it is biflat and has a one-sided b.a.i. On the other hand, if $A$ is biprojective
and has a one-sided identity, then it is superbiprojective \cite{Sel_super}.

\section{Biprojectivity}
\begin{lemma}
\label{lemma:split_ext}
Let 
\begin{equation}
\label{ext1}
0 \rar I \rar \A \xra{\sigma} A\rar 0
\end{equation}
be an extension of Banach algebras
such that $\A I=0$. Assume there exists an antihomomorphism
$\beta\colon A\to A$ and a linear continuous map $\alpha\colon \A\to\A$
such that $\beta^2=\id_A$ and $\beta\sigma=\sigma\alpha$.
Suppose also that $A$ is essential and projective as a
right $\A$-module via $\sigma$.
Then \eqref{ext1} is admissible.
\end{lemma}
\begin{proof}
Denote by $\pi\colon A\Ptens\A\to A$ the right action of $\A$ on $A$ determined
by $\sigma$, i.e., $\pi(a\otimes u)=a\sigma(u)$.
Since $A$ is essential and projective in $\rmod\A$, there exists
a right $\A$-module morphism $\rho\colon A\to A\Ptens\A$ such that
$\pi\rho=\id_A$ \cite[IV.1]{X1}. Next observe that the condition $\A I=0$ implies
that $\A$ has a natural structure of right Banach $A$-module. Indeed,
the product map $\A\times\A\to\A$ vanishes on $\A\times I$ and hence
determines a right action $\A\times A=\A\times (\A/I)\to\A$ by the rule
$(u,a)\mapsto u\cdot a=uv$ where $v\in\sigma^{-1}(a)$.
The corresponding linear map $\A\Ptens A\to\A,\; u\otimes a\mapsto u\cdot a$
will be denoted by $\phi$.

Define $\tilde\varkappa\colon A\to\A$ as the composition
\begin{equation*}
A \xra{\rho} A\Ptens\A \xra{\tau} \A\Ptens A \xra{\alpha\otimes\beta}
\A\Ptens A \xra{\phi} \A
\end{equation*}
where $\tau$ stands for the flip $a\otimes u\mapsto u\otimes a$.

Let us compute $\sigma\tilde\varkappa$.
For every $u\in\A,\; a\in A$, and $v\in\sigma^{-1}(a)$ we have
\begin{equation*}
(\sigma\phi)(u\otimes a)=\sigma(uv)=\sigma(u)\sigma(v)=\sigma(u)a\, .
\end{equation*}
In other words, $\sigma\phi=\pi^{\op}$, where $\pi^{\op}\colon \A\Ptens A\to A$
is the left action of $\A$ on $A$ determined by $\sigma$.   
Next,
\begin{multline*}
\bigl(\pi^{\op}\circ(\alpha\otimes\beta)\bigr)(u\otimes a)=
\pi^{\op}\bigl(\alpha(u)\otimes\beta(a)\bigr)=
\sigma\bigl(\alpha(u)\bigr)\beta(a)\\
=\beta\bigl(\sigma(u)\bigr)\beta(a)=
\beta\bigl(a\sigma(u)\bigr)=
(\beta\pi)(a\otimes u)=
(\beta\pi\tau)(u\otimes a)\, .
\end{multline*}
Hence $\pi^{\op}\circ(\alpha\otimes\beta)=\beta\pi\tau$. Finally,
\begin{equation*}
\sigma\tilde\varkappa=
\pi^{\op}\circ (\alpha\otimes\beta)\tau\rho=
\beta\pi\tau\tau\rho=\beta\pi\rho=\beta\, .
\end{equation*}
Since $\beta^2=\id_A$, we conclude that the map $\varkappa=\tilde\varkappa\beta$
satisfies $\sigma\varkappa=\id_A$. Therefore \eqref{ext1} is admissible.
\end{proof}

Recall (see, e.g., \cite{DU,WBJhnsn}) that a Banach space $E$ is said to
have the {\em Radon-Nikod\'ym property} (RNP for short) if
for each finite measure space $(X,\mu)$ every $\mu$-continuous
$E$-valued measure of finite variation is differentiable w.r.t. $\mu$. 
For our purposes, the following properties of the RNP will be important.
\begin{itemize}
\item[(a)] The RNP is inherited by closed subspaces.
\item[(b)] If $E$ and $F$ are reflexive Banach spaces one of which
has the approximation property, then the space
$\mathscr N(E^*,F)$ of nuclear operators from $E^*$ to $F$ has the RNP
\cite{Diest_Fair}. In partcular, $\nN(L^p(X,\mu))$ has the RNP for every
measure space $(X,\mu)$. 
\item[(c)] If $(X,\mu)$ is a measure space, then
$L^1(X,\mu)$ has the RNP if and only if $\mu$ is purely atomic
\cite[III.1]{DU}.
\end{itemize}
 
\begin{lemma}
\label{lemma:neuf_ext}
Let $G$ be a locally compact group and let $1<p<\infty$.
Then the following conditions are equivalent:
\begin{itemize}
\item[(i)] Extension \eqref{neuf_ext} is admissible;
\item[(ii)] Extension \eqref{neuf_ext} splits;
\item[(iii)] $G$ is discrete.
\end{itemize}
\end{lemma}
\begin{proof}
(ii)$\Rightarrow$(i): obvious.

(i)$\Rightarrow$(iii).
If $G$ is nondiscrete, then $L^1(G)$ does not have the RNP
(see (c) above).
Since $\nN^p(G)=\nN(L^p(G))$ has the RNP (see (b) above), we conclude
that $L^1(G)$ is not isomorphic to a subspace of $\nN^p(G)$.
Hence extension \eqref{neuf_ext} is not admissible.

(iii)$\Rightarrow$(ii).  
If $G$ is discrete, then the map
\begin{equation*}
\varkappa\colon\ell^1(G)\to\nN^p(G),\qquad f\mapsto M_{\check f}
\end{equation*}
is a continuous right inverse to $\sigma\colon\nN^p(G)\to\ell^1(G)$
(see \cite{Neuf_PhD}, Satz 5.3.7).
It is easy to check that $\varkappa$ is an algebra
homomorphism. To see this, for each $t\in G$ let
$\delta_t$ denote the function which equals $1$ at $t$ and $0$ elsewhere.
We claim that $M_{\delta_s}*M_{\delta_t}=M_{\delta_{ts}}$ for each $s,t\in G$.
Indeed, for each $T\in\B(\ell^p(G))$ we have
\begin{multline*}
\langle T, M_{\delta_s}*M_{\delta_t}\rangle
=\langle T\odot M_{\delta_s},\delta_t\otimes\delta_t\rangle
=\langle (T\odot M_{\delta_s})(\delta_t),\delta_t\rangle
=\langle \delta_s\otimes\delta_s, L_t T L_{t^{-1}} \rangle\\
=\langle T L_{t^{-1}}(\delta_s),L_{t^{-1}}(\delta_s)\rangle
=\langle T(\delta_{ts}),\delta_{ts}\rangle
=\langle T, \delta_{ts}\otimes\delta_{ts} \rangle
=\langle T,M_{\delta_{ts}}\rangle\, .
\end{multline*}
Therefore $M_{\delta_s}*M_{\delta_t}=M_{\delta_{ts}}$ for each $s,t\in G$, and so
$$
\varkappa(\delta_s*\delta_t)=
\varkappa(\delta_{st})=M_{\delta_{t^{-1}s^{-1}}}
=M_{\delta_{s^{-1}}}*M_{\delta_{t^{-1}}}
=\varkappa(\delta_s)*\varkappa(\delta_t)\, .
$$
Thus we see that $\varkappa$ is an algebra homomorphism.
Next,
$$
\sigma\varkappa(\delta_s)=\sigma(\delta_{s^{-1}}\otimes\delta_{s^{-1}})=\delta_s,
$$
i.e., $\sigma\varkappa=\id_{\ell^1(G)}$ (see also \cite{Neuf_PhD}, Satz 5.3.7).
Therefore extension \eqref{neuf_ext} splits.
\end{proof}

\begin{lemma}
\label{lemma:proj_quot}
Let \eqref{ext1} be an extension of Banach algebras such that $\A I=0$.
Assume that \eqref{ext1} splits and $A$ has a left identity. Then
$A$ is projective in $\rmod\A$.
\end{lemma}
\begin{proof}
As in Lemma \ref{lemma:split_ext}, denote by  
$\pi\colon A\Ptens\A\to A$ the right action of $\A$ on $A$ given by
$\pi(a\otimes u)=a\sigma(u)$. Let $\varkappa\colon A\to\A$ be a homomorphism
such that $\sigma\varkappa=\id_A$. Take a left identity $e$ of $A$ and
define the map $\rho\colon A\to A\Ptens\A$
by the rule $\rho(a)=e\otimes\varkappa(a)$.
Evidently, $\pi\rho=\id_A$.
Since $\A I=0$, we see that $vu=v\varkappa(\sigma(u))$
for each $u,v\in\A$. Hence
$$
\rho(a\cdot u)=\rho(a\sigma(u))=e\otimes\varkappa(a\sigma(u))=
e\otimes\varkappa(a)\varkappa(\sigma(u))=
e\otimes\varkappa(a)u=\rho(a)\cdot u
$$
for each $a\in A$ and each $u\in\A$. This means that $\rho$ is a morphism
in $\rmod\A$.
Since $\pi\rho=\id_A$, we conclude that $A$ is projective in $\rmod\A$.
\end{proof}

\begin{theorem}
\label{thm:discr}
Let $G$ be a locally compact group and let $1<p<\infty$.
Then the following conditions are equivalent:
\begin{itemize}
\item[(i)]
$L^1(G)$ is projective in $\rmod\nN^p(G)$;
\item[(ii)]
$G$ is discrete.
\end{itemize}
\end{theorem}
\begin{proof}
(i)$\Rightarrow$(ii).
For each $1\le r\le\infty$ define the map $\alpha_r\colon L^r(G)\to L^r(G)$ by
$\alpha_r(f)(t)=\Delta(t^{-1})^{1/r} f(t^{-1})$. It is easy to check that
$\alpha_r$ is an isometry and $\alpha_r^2=\alpha_r$.
Now let $\beta=\alpha_1\colon L^1(G)\to L^1(G)$
and $\alpha=\alpha_p\otimes\alpha_q\colon\nN^p(G)\to\nN^p(G)$.
Evidently, $\beta$ is an antihomomorphism.
For each $f\in L^p(G)$ and each $g\in L^q(G)$ we have
\begin{multline*}
\sigma(\alpha(f\otimes g))(t)
=\sigma(\alpha_p(f)\otimes\alpha_q(g))(t)
=\Delta(t^{-1}) \alpha_p(f)(t^{-1}) \alpha_q(g)(t^{-1})\\
=\Delta(t)^{-1} \Delta(t)^{1/p} f(t) \Delta(t)^{1/q} g(t)
=f(t) g(t)\, .
\end{multline*}
On the other hand,
$$
\beta(\sigma(f\otimes g))(t)
=\Delta(t^{-1}) \sigma(f\otimes g)(t^{-1})
=\Delta(t)^{-1}\Delta(t)f(t) g(t)=f(t) g(t)\, .
$$
Hence $\beta\sigma=\sigma\alpha$. Finally, $L^1(G)$ is an essential
$\nN^p(G)$-module since $\sigma$ is surjective and $L^1(G)$ has a b.a.i.
By Lemma \ref{lemma:split_ext}, extension \eqref{neuf_ext} is
admissible. Now Lemma \ref{lemma:neuf_ext} shows that $G$ is discrete. 

(ii)$\Rightarrow$(i).
If $G$ is discrete, then extension \eqref{neuf_ext} splits
by Lemma \ref{lemma:neuf_ext}. Since $L^1(G)$ is unital in this case,
Lemma \ref{lemma:proj_quot} implies
that $L^1(G)$ is projective in $\rmod\nN^p(G)$.
\end{proof} 

\begin{lemma}
\label{lemma:conv_comp}
Let $G$ be a compact group and let $1$ denote the function that is
identically $1$ on $G$. Then 
$(1\otimes 1)*a=\Tr a\cdot 1\otimes 1$ for each $a\in\nN^p(G)$.
\end{lemma}
\begin{proof}
For each $T\in\B(L^p(G))$ we have $T\odot(1\otimes 1)=M_h$, where
\begin{multline*}
h(t)=\langle 1\otimes 1, L_t T L_{t^{-1}}\rangle
=\langle L_t T L_{t^{-1}}(1), 1\rangle
=\langle T L_{t^{-1}}(1), L_{t^{-1}}(1)\rangle
=\langle T(1), 1\rangle
\end{multline*}
for each $t\in G$. Hence $M_h=\langle T(1), 1\rangle \id_{L^p(G)}$,
and for each $a\in\nN^p(G)$ we have
\begin{multline*}
\langle T, (1\otimes 1)*a\rangle
=\langle T\odot(1\otimes 1),a\rangle
=\langle T(1), 1\rangle \langle \id_{L^p(G)}, a\rangle\\
=\langle T(1), 1\rangle\Tr a =
\langle T, \Tr a\cdot 1\otimes 1\rangle\, ,
\end{multline*}
as required.
\end{proof}

\begin{theorem}
\label{thm:comp}
Let $G$ be a locally compact group and let $1<p<\infty$.
Then the following conditions are equivalent:
\begin{itemize}
\item[(i)]
$\CC_{\Tr}$ is projective in $\rmod\nN^p(G)$;
\item[(ii)]
$G$ is compact.
\end{itemize}
\end{theorem}
\begin{proof}
(i)$\Rightarrow$(ii). If $\CC_{\Tr}$ is projective in $\rmod\nN^p(G)$, then
$\CC_\varepsilon$ is projective in $\rmod L^1(G)$ (see Remark \ref{rem:pullback}).
This means exactly that $G$ is compact (see \cite{X_flat} or \cite[IV.5]{X1}).

(ii)$\Rightarrow$(i). First observe that the canonical morphism
$$
\pi\colon\CC_{\Tr}\Ptens\nN^p(G)\to\CC_{\Tr},\quad
\lambda\otimes u\mapsto\lambda\cdot u,
$$
is identified with $\Tr\colon\nN^p(G)\to\CC_{\Tr}$.
Define $\rho\colon\CC_{\Tr}\to\nN^p(G)$ by
$\rho(\lambda)=\lambda\cdot 1\otimes 1$. It is clear that $\pi\rho=\id_{\CC}$,
and Lemma \ref{lemma:conv_comp} implies that $\rho$ is a right
$\nN^p(G)$-module morphism. Hence $\CC_{\Tr}$ is projective
in $\rmod\nN^p(G)$.
\end{proof}

\begin{theorem}
\label{thm:bipr}
Let $G$ be a locally compact group and let $1< p <\infty$.
Then the following conditions are equivalent:
\begin{itemize}
\item[(i)]
$\nN^p(G)$ is biprojective;
\item[(ii)]
$\nN^p(G)$ is superbiprojective;
\item[(iii)]
$G$ is finite.
\end{itemize}
\end{theorem}
\begin{proof}
(ii)$\Rightarrow$(i): obvious.

\noindent
(iii)$\Rightarrow$(ii). If $G$ is finite, then $L^1(G)=\CC[G]$ is contractible.
Furthermore, all the algebras in the extension \eqref{neuf_ext} are
finite-dimensional. Hence \eqref{neuf_ext} splits, and so $\nN^p(G)$
is isomorphic to the semidirect product $L^1(G)\oplus I$.   
Since $L^1(G)\cdot I=0$, we conclude that $\nN^p(G)$ is superbiprojective
by \cite{Sel_super} (see also Lemma \ref{lemma:bipr_ext} below).

\noindent
(i)$\Rightarrow$(iii). 
If $\nN^p(G)$ is biprojective, then both $L^1(G)$ and $\CC_{\Tr}$
are projective in $\rmod\nN^p(G)$ since they are essential and
$\nN^p(G)$ has a right b.a.i. (see \cite{Sel_bipr} or \cite[7.1.60]{X2}).
Now it remains to apply Theorems \ref{thm:discr} and \ref{thm:comp}.
\end{proof}
 
\section{Biflatness}

\begin{lemma}
\label{lemma:bifl_ext}
Let
\begin{equation}
\label{ext2}
0 \rar I \rar \A \xra{\sigma} A\rar 0
\end{equation}
be a weakly admissible extension of Banach algebras
such that $\A I=0$ and $\A^2=\A$. Then the following conditions
are equivalent:
\begin{itemize}
\item[(i)]
$\A$ is biflat and has a right b.a.i.;
\item[(ii)]
$A$ is biflat and has a right b.a.i.
\end{itemize}
\end{lemma}
\begin{proof}
(i)$\Rightarrow$(ii). Since $\A$ has a right b.a.i., we have $I\A=I$.
Hence $A=\A/(I\A)$ is biflat by \cite{Sel_wdb}. It is also clear that
$A$ has a right b.a.i.

(ii)$\Rightarrow$(i). 
As in Lemma \ref{lemma:split_ext}, the condition $\A I=0$ implies that
$\A$ is a right Banach $A$-module in a natural way.
First we note that $\A$ has a right b.a.i. Indeed, 
let $\{ e_\nu\}$ be a bounded net in $\A$ such that
$\{\sigma(e_\nu)\}$ is a right b.a.i. in $A$. Then for every $u,v\in\A$ we have
$$
\| uve_\nu-uv\|
=\| u\cdot\sigma(ve_\nu)-u\cdot\sigma(v) \|
\le \| u\| \|\sigma(v)\sigma(e_\nu)-\sigma(v)\| \to 0\, ,
$$
because $\{ \sigma(e_\nu)\}$ is a right b.a.i. in $A$. Since $\{ e_\nu\}$
is bounded, and since $\A^2=\A$, we conclude that $\{ e_\nu\}$ is a right
b.a.i. in $\A$. 

Let $\pi_A\colon A\Ptens A\to A$ be the product map. Since $A$ is biflat,
the dual map $\pi_A^*$ is a coretraction in $A\bimod A$.
Applying the functor $\h_A(\A,?)$, we see that
$\h_A(\A,\pi_A^*)\colon\h_A(\A,A^*)\to\h_A(\A,(A\Ptens A)^*)$
is a coretraction in $A\bimod\A$. Using the adjoint associativity isomorphisms
\cite[II.5]{X1}, we can identify the latter map with
$(\id_\A\otimes\pi_A)^*\colon (\A\ptens{A} A)^*\to (\A\ptens{A} (A\Ptens A))^*$.
Therefore $(\id_\A\otimes\pi_A)^*$ is a coretraction in $A\bimod\A$.

Now let $\pi_{\A,A}\colon \A\Ptens A\to\A,\; u\otimes a\mapsto u\cdot a$
denote the right action on $A$ on $\A$.
Evidently, $\pi_{\A,A}$ induces the morphism
$\varkappa\colon \A\ptens{A} A\to\A,\; u\otimes a\mapsto u\cdot a$.
We have the following commutative diagram in $\A\bimod A$:
\begin{equation}
\label{pi_1}
\begin{CD}
\A\ptens{A} A\Ptens A @> {\id_\A\otimes\pi_A} >> \A\ptens{A} A\\
@V {\varkappa\otimes\id_A} VV @VV {\varkappa} V\\
\A\Ptens A @> {\pi_{\A,A}} >> \A
\end{CD}
\end{equation}
The vertical arrows in the diagram are isomorphisms in 
$\A\bimod A$, because $\A^2=\A$ and
$A$ has a right b.a.i. \cite[II.3]{X1}. 
Consider now the dual diagram. We already know that
$(\id_\A\otimes\pi_A)^*$ is a coretraction in $A\bimod\A$.
Hence so is $\pi_{\A,A}^*$.

Since \eqref{ext2} is weakly admissible, 
$\sigma^*$ is an admissible monomorphism in $A\lmod$.
On the other hand, since $A$ is biflat, it follows that $A$ is flat in $\rmod A$,
and hence $A^*$ is injective
in $A\lmod$ (see, e.g., \cite[VII.1]{X1}).
Therefore $\sigma^*$ is a coretraction in $A\lmod$, and so
$\B(\A,\sigma^*)\colon\B(\A,A^*)\to\B(\A,\A^*)$ is a coretraction
in $A\bimod\A$.
Identifying $\B(\A,A^*)$ with $(\A\Ptens A)^*$ and $\B(\A,\A^*)$
with $(\A\Ptens\A)^*$, we conclude that
$(\id_\A\otimes\sigma)^*\colon (\A\Ptens A)^*\to (\A\Ptens\A)^* $
is a coretraction in $A\bimod\A$.

It is easy to see that the product map $\pi_\A\colon \A\Ptens\A\to\A$
decomposes as $\pi_\A=\pi_{\A.A}\circ (\id_\A\otimes\sigma)$.
Hence $\pi_\A^*=(\id_\A\otimes\sigma)^*\circ \pi_{\A.A}^*$.
But we already know that both $(\id_\A\otimes\sigma)^*$ and $\pi_{\A.A}^*$
are coretractions in $A\bimod\A$.
Hence so is $\pi_\A^*$. 

To complete the proof, it remains to show that $\pi_\A^*$ is actually
a coretraction in $\A\bimod\A$. To this end, note that
for each $f\in\A^*$ and each $u,v\in\A$ we have
$$
\langle u\cdot f,v\rangle =
\langle f,vu\rangle=
\langle f,v\cdot\sigma(u)\rangle=
\langle\sigma(u)\cdot f,v\rangle\, ,
$$
i.e., $u\cdot f=\sigma(u)\cdot f$ for each $f\in\A^*$ and each $u\in\A$.
Similarly, $u\cdot g=\sigma(u)\cdot g$ for each $g\in(\A\Ptens\A)^*$.
This implies that every left $A$-module morphism between $(\A\Ptens\A)^*$
and $\A^*$ is in fact a left $\A$-module morphism. In particular,
a left inverse of $\pi_\A^*$ in $A\bimod\A$ is a morphism in $\A\bimod\A$.
Therefore $\pi_\A^*$ is a coretraction in $\A\bimod\A$. 
This completes the proof.
\end{proof}

Lemma \ref{lemma:bifl_ext} has the following ``predual'' counterpart,
which is a slight generalization of an unpublished result of Selivanov \cite{Sel_super}.

\begin{lemma}
\label{lemma:bipr_ext} 
Let
\begin{equation}
\label{ext3}
0 \rar I \rar \A \xra{\sigma} A\rar 0
\end{equation}
be an admissible extension of Banach algebras
such that $\A I=0$ and $\A^2=\A$. Then the following conditions
are equivalent:
\begin{itemize}
\item[(i)]
$\A$ is biprojective and has a right b.a.i. (resp. a right identity);
\item[(ii)]
$A$ is biprojective and has a right b.a.i. (resp. a right identity).
\end{itemize}
As a corollary, $\A$ is superbiprojective provided $A$ is contractible.
\end{lemma}
\begin{proof}
The proof is similar to that of Lemma \ref{lemma:bifl_ext}.

(i)$\Rightarrow$(ii). Since $\A$ has a right b.a.i., we have $I\A=I$.
Hence $A=\A/(I\A)$ is biprojective by \cite{Sel_bipr}. Evidently,
if $\A$ has a right b.a.i. (resp. a right identity), then so does $A$.

(ii)$\Rightarrow$(i). 
As in Lemma \ref{lemma:split_ext}, the condition $\A I=0$ implies that
$\A$ is a right Banach $A$-module in a natural way.
Arguing as in the proof of Lemma \ref{lemma:bifl_ext}, we conclude that
any bounded preimage of a right b.a.i. (resp. of a right identity) in $A$ is a 
right b.a.i. (resp. a right identity) in $\A$.

Since $A$ is biprojective, the product map
$\pi_A\colon A\Ptens A\to A$ is a retraction in $A\bimod A$.
Hence $\id_\A\otimes\pi_A\colon\A\ptens{A}(A\Ptens A)\to\A\ptens{A} A$
is a retraction in $\A\bimod A$.
As in Lemma \ref{lemma:bifl_ext}, we can identify the latter morphism
with $\pi_{\A,A}\colon \A\Ptens A\to\A$ (see diagram \eqref{pi_1}).
Therefore $\pi_{\A,A}$ is also a retraction in $\A\bimod A$.

Since \eqref{ext3} is admissible, 
$\sigma$ is an admissible epimorphism in $\rmod A$.
On the other hand, since $A$ is biprojective, it follows that $A$ is
projective in $\rmod A$. Hence $\sigma$ is a retraction in $\rmod A$,
and so $\id_\A\otimes\sigma\colon \A\Ptens\A\to \A\Ptens A$
is a retraction in $\A\bimod A$.

It is easy to see that the product map $\pi_\A\colon \A\Ptens\A\to\A$
decomposes as $\pi_\A=\pi_{\A.A}\circ (\id_\A\otimes\sigma)$.
But we already know that both $\id_\A\otimes\sigma$ and $\pi_{\A.A}$
are retractions in $\A\bimod A$. Hence so is $\pi_\A$. 

Finally, the condition $\A I=0$ implies that $uv=u\cdot\sigma(v)$
and $w\cdot v=w\cdot\sigma(v)$ for each $u,v\in\A$ and each $w\in\A\Ptens\A$.
Hence every right $A$-module morphism between $\A\Ptens\A$
and $\A$ is in fact a right $\A$-module morphism. In particular,
a right inverse of $\pi_\A$ in $\A\bimod A$ is a morphism in $\A\bimod\A$.
Therefore $\pi_\A$ is a retraction in $\A\bimod\A$, i.e., $\A$ is biprojective. 
This completes the proof.
\end{proof}

\begin{theorem}
\label{thm:biflat}
Let $G$ be a locally compact group and let $1< p <\infty$.
Then the following conditions are equivalent:
\begin{itemize}
\item[(i)]
$\nN^p(G)$ is biflat;
\item[(ii)]
$\nN^p(G)$ is superbiflat;
\item[(iii)]
$\CC_{\Tr}$ is flat in $\rmod\nN^p(G)$;
\item[(iv)]
$G$ is amenable.
\end{itemize}
\end{theorem}
\begin{proof}
(i)$\iff$(ii): clear, because $\nN^p(G)$ has a right b.a.i.

(i)$\Rightarrow$(iii). If $\nN^p(G)$ is biflat, then $\CC_{\Tr}$ is flat
in $\rmod\nN^p(G)$ since $\CC_{\Tr}$ is essential and
$\nN^p(G)$ has a right b.a.i. (see \cite[7.1.60]{X2}).

(iii)$\Rightarrow$(iv). If $\CC_{\Tr}$ is flat in $\rmod\nN^p(G)$, then
$\CC_\varepsilon$ is flat in $\rmod L^1(G)$ (see Remark \ref{rem:pullback}).
By \cite{X_flat} (see also \cite[VII.2]{X1}),
this happens if and only if $G$ is amenable.

(iv)$\Rightarrow$(i). Recall that $G$ is amenable if and only 
if $L^1(G)$ is amenable (see, e.g., \cite{Jhnsn_CBA,X_flat,X1}).
Since $L^\infty(G)=L^1(G)^*$ is an injective Banach space (see, e.g.,
\cite{Wojt}), we see that extension \eqref{neuf_ext} is weakly
admissible. Therefore $\nN^p(G)$ is biflat by Lemma \ref{lemma:bifl_ext}.
\end{proof}  

\begin{remark}
M. Neufang has kindly informed the author that he has also proved the equivalence
of conditions (i), (ii) and (iv) of Theorem \ref{thm:biflat}.
\end{remark} 

\noindent\textbf{Acknowledgments. }The author is grateful to
A.~Ya.~Helemskii and Yu.~V.~Selivanov for valuable discussions.

%\bigskip
\vfill
\begin{flushleft}
\scshape\small
Department of Differential Equations and Functional Analysis\\
Faculty of Physics, Mathematics, and Sciences\\
Peoples' Friendship University of Russia\\
Mikluho-Maklaya 6\\
117198 Moscow\\
RUSSIA

\medskip
{\itshape Address for correspondence:}\\

\medskip\upshape
Krupskoi 8--3--89\\
Moscow 119311\\
Russia

\medskip
{\itshape E-mail:} {\ttfamily pirkosha@sci.pfu.edu.ru, pirkosha@online.ru}
\end{flushleft}

\begin{thebibliography}{99}
\bibitem{Dales}
Dales, H. G.
{\em Banach Algebras and Automatic Continuity},
Clarendon Press, Oxford, 2000.
\bibitem{Diest_Fair}
Diestel, J. and Faires, B.
{\em On vector measures},
Trans. Amer. Math. Soc. \textbf{198} (1974), 253--271.
\bibitem{DU}
Diestel, J. and Uhl, J. J., Jr.
{\em Vector measures},
Amer. Math. Soc., Providence, 1977.
\bibitem{HRII}
Hewitt, E. and Ross, K. A.
{\em Abstract harmonic analysis, Vol.~II},
Springer, 1970.
\bibitem{X_calc_est}
Helemskii, A. Ya.
{\em On a method for calculating and estimating
the global homological dimension of Banach algebras},
Mat. Sbornik \textbf{87} (129) (1972), 122--135 (Russian);
English transl.: Math. USSR Sb. \textbf{16} (1972), 125--138.
\bibitem{X_flat}
Helemskii, A. Ya.
{\em Flat Banach modules and amenable algebras},
Trudy MMO \textbf{47} (1984), 179--218 (Russian);
English transl.: Trans. Moscow Math. Soc. \textbf{47} (1985), 199--244.
\bibitem{X1}
Helemskii, A. Ya.
{\itshape The Homology of Banach and Topological Algebras},
Moscow University Press, 1986 (Russian);
English transl.: Kluwer Academic Publishers, Dordrecht, 1989.
\bibitem{X2}
Helemskii, A. Ya.
{\itshape Banach and Polynormed Algebras: General Theory,
Representations, Homology},
Nauka, Moscow, 1989 (Russian); English transl.:
Oxford University Press, 1993.
\bibitem{Jhnsn_CBA}
Johnson, B. E.
{\em Cohomology in Banach algebras},
Mem. Amer. Math. Soc. \textbf{127} (1972).
\bibitem{WBJhnsn}
Johnson, W. B. and Lindenstrauss, J.
{\em Basic concepts in the geometry of Banach spaces},
Handbook of the geometry of Banach spaces, Vol. I, 1--84, 
North-Holland, Amsterdam, 2001. 
\bibitem{Neuf_PhD}
Neufang, M.
{\em Abstrakte Harmonische Analyse und
Modulhomomorphismen \"uber von Neumann-Algebren},
Dissertation zur Erlangung des Grades des Doktors
des Naturwissenschaften, Saarbr\"ucken, 2000.
\bibitem{Runde}
Runde, V.
{\em Lectures on Amenability},
Lecture Notes in Math. 1774, Springer, 2002.
\bibitem{Sel_bipr}
Selivanov, Yu. V.
{\em Biprojective Banach algebras},
Izv. Akad. Nauk SSSR ser. mat. \textbf{43} (1979), 1159--1174;
English transl.: Math. USSR Izvestija \textbf{15} (1980), 387--399.
\bibitem{Sel_wdb}
Selivanov, Yu. V.
{\em Weak homological bidimension and its values in
the class of biflat Banach algebras},
Extracta Math. \textbf{11} (1996), 348--365.
\bibitem{Sel_super}
Selivanov, Yu. V.
{\em Superbiprojective and superbiflat Banach algebras},
Unpublished manuscript, Odense, 2001.
\bibitem{Wendel}
Wendel, J. G.
{\em Left centralizers and isomorphisms of group algebras},
Pacific J. Math. \textbf{2} (1952), 251--261.
\bibitem{Wojt}
Wojtaszczyk, P. 
{\em Banach spaces for analysts}, 
Cambridge University Press, Cambridge, 1991.  
\end{thebibliography}
\end{document}